\documentclass{amsart}

\usepackage{amsmath} 
\usepackage{amssymb}

\newtheorem{theorem}{Theorem}[section] 
\newtheorem{claim}{Claim}[theorem]
\newtheorem{lemma}[theorem]{Lemma} 
\newtheorem{proposition}[theorem]{Proposition} 
\newtheorem{observation}[theorem]{Observation} 
\newtheorem{corollary}[theorem]{Corollary} 

\theoremstyle{definition}
\newtheorem{definition}[theorem]{Definition}

\newtheorem{problem}[theorem]{Problem}
\newtheorem{aim}[theorem]{Aim}

\theoremstyle{remark}
\newtheorem{remark}[theorem]{Remark}

\numberwithin{equation}{section}
\newcommand{\forces}{\Vdash}

\newcommand{\bV}{{\bf V}} 
\newcommand{\lesdot}{\mathrel{\mathord{<}\!\!\raise 
0.8 pt\hbox{$\scriptstyle\circ$}}}

\newcommand{\lh}{{\rm lh}\/}
\newcommand{\rest}{{\restriction}}

\newcommand{\cA}{{\mathcal A}}

\newcommand{\cF}{{\mathcal F}}

\newcommand{\bbP}{{\mathbb P}}
\newcommand{\cP}{{\mathcal P}}

\newcommand{\cT}{{\mathcal T}}
\newcommand{\cU}{{\mathcal U}}

\newcommand{\cl}{{\rm cl}\/} 
\newcommand{\fil}{{\rm fil}\/} 
 
\newcommand{\ulf}{{\rm ulf}\/} 
\newcommand{\uuf}{{\rm uuf}\/} 
\newcommand{\otp}{{\rm otp}\/} 
\newcommand{\cFl}{{\mathcal F}_\lambda}
\newcommand{\st}{{\bf st}} 
\newcommand{\cHl}{{\mathcal H}(\lambda^{++})}
\newcommand{\vare}{\varepsilon}
\newcommand{\gdl}{{\mathfrak d}_\lambda} 
\newcommand{\gcl}{{\mathfrak d}_\cl(\lambda)} 
\newcommand{\bql}{{{\mathbb Q}^1_\lambda}}
\newcommand{\bqz}{{{\mathbb Q}^0_\lambda}}

\newcommand{\UE}{{\rm UE}}
\newcommand{\UF}{{\rm UF}}
\newcommand{\Mkkk}{{\bf M}^\lambda_{\lambda,\lambda}}

\newcount\skewfactor
\def\mathunderaccent#1#2 {\let\theaccent#1\skewfactor#2
\mathpalette\putaccentunder}
\def\putaccentunder#1#2{\oalign{$#1#2$\crcr\hidewidth
\vbox to.2ex{\hbox{$#1\skew\skewfactor\theaccent{}$}\vss}\hidewidth}}
\def\name{\mathunderaccent\tilde-3 }

\begin{document}

\title{The combinatorics of reasonable ultrafilters}

\author{Saharon Shelah}
\address{Einstein Institute of Mathematics\\
Edmond J. Safra Campus, Givat Ram\\
The Hebrew University of Jerusalem\\
Jerusalem, 91904, Israel\\
 and  Department of Mathematics\\
 Rutgers University\\
 New Brunswick, NJ 08854, USA}
\email{shelah@math.huji.ac.il}
\urladdr{http://shelah.logic.at}
\thanks{The author acknowledges support from the United States-Israel
Binational Science Foundation (Grant no. 2002323). Publication 830.}     

\subjclass{Primary 03E05; Secondary: 03E20}
\date{September 2005}

\begin{abstract}
We are interested in generalizing part of the theory of ultrafilters
on $\omega$ to larger cardinals. Here we set the scene for further
investigations introducing properties of ultrafilters in strong sense dual
to being normal. 
\end{abstract}

\maketitle

\section{Introduction}
Questions concerning ultrafilters on $\omega$ have occurred to be very 
stimulating for research in several subareas of Set Theory and Topology. We
hope that this success story could be repeated for ultrafilters on
uncountable regular cardinals $\lambda$, particularly if $\lambda$ is
strongly inaccessible. Our aim in the present paper is to introduce new 
properties of ultrafilters and argue that these properties could play the
stimulating role that was once played by {\em $P$--points on $\omega$}. 

In a long run, we plan to find generalizations of the following results:
\begin{enumerate}
\item[(a)]  Consistently, some ultrafilters on $\omega$ are generated by
  $<2^{\aleph_0}$ many sets.
\item[(b)]  $P$-points are preserved by some forcing notions (see, e.g., 
\cite[V]{Sh:f}, \cite{RoSh:470}).
\item[(c)] Consistently, there is no $P$--point.
\item[(d)]  For a function $f:\omega\longrightarrow\omega$ and ultrafilter
 $D$ on $\omega$, let 
\[D/f\stackrel{\rm def}{=}\{A\subseteq\omega:f^{-1}(A)\in D\};\]
it is an ultrafilter on $\omega$ (of course, we are interested in the cases 
when $D$ and $D/f$ are uniform, which in this case is the same as
non-principal). By Blass and Shelah \cite{BsSh:242}, consistently for any
two non-principal ultrafilters $D_1,D_2$ on $\omega$ there are finite-to-one
non-decreasing functions $f_1,f_2:\omega\longrightarrow\omega$ such that
$D_1/f_1 = D_2/f_2$. 
\item[(f)]  For a significant family of forcing notions built according to
the scheme of creatures of \cite{RoSh:470} we can consider an appropriate
filter, i.e., if $\langle p_\alpha:\alpha <\omega_1 \rangle$ is
$\le^*$-increasing it may define an ultrafilter which is not necessarily
generated by $\aleph_1$-sets, so we may ask on this. 
\end{enumerate}

There is much work on normal ultrafilters, the parallel on
$\omega$ are Ramsey ultrafilters. Now, every Ramsey ultrafilter on $\omega$
is a $P$-point but there are $P$-points of very different characters, e.g., 
$P$-point with no Ramsey ultrafilter below.   Gitik \cite{Gi81} has
investigated generalizations of $P$-points for normal ultrafilters. But this 
paper goes in a different direction (which up to recently I have not 
considered to be fruitful) and we restrict our attention to ultrafilters
which are very non-normal --- {\em the weakly reasonable ultrafilters}. What
is a weakly reasonable ultrafilter on $\lambda$?  It is a uniform
ultrafilter on a regular cardinal $\lambda$ which does not contain some club
of $\lambda$ and such that this property is preserved if we divide it by a
non-decreasing $f:\lambda \longrightarrow \lambda$ with unbounded range (see
Definition \ref{1.5} below).  

We also want that our ultrafilters generalize $P$--points on $\omega$ and in
the second section we introduce {\em reasonable\/} and {\em very reasonable}
ultrafilters. The property defining $P$--points is that countable families
of sets from the ultrafilter have pseudo-intersections in the
ultrafilter. We modify this property so that we involve some description of
how the considered ultrafilter is generated, and we postulate that {\em the
generating systems\/} are suitably directed. This is a replacement for the 
existence of pseudo-intersections and it is the essence 
of Definition \ref{1.3}(4,5). The third section shows that the number of
generating systems (of our type) for somewhat reasonable ultrafilters cannot
be too small. We conclude the paper with a section listing open problems and
describing further research. 
\medskip

\noindent{\bf Notation:}\quad Our notation is rather standard and compatible
with that of classical textbooks (like Jech \cite{J}). In forcing we keep
the older convention that {\em a stronger condition is the larger
one}. (However, in the present paper we use forcing notions only for
combinatorial constructions and almost every mention of forcing just means
that we a dealing with a transitive reflexive relation
$\bbP=(\bbP,\leq_\bbP)$.) 

\begin{enumerate}
\item Ordinal numbers will be denoted be the lower case initial letters of
the Greek alphabet ($\alpha,\beta,\gamma,\delta\ldots$) and also by $i,j$
(with possible sub- and superscripts). 
\item Cardinal numbers will be called $\kappa,\lambda,\mu$; $\lambda$ will be
always assumed to be a regular uncountable cardinal (we may forget to
mention it). 
\item $D,\cU$ will denote filters on $\lambda$, $G,G^*,G^*_\ell$ will be
  subsets of specific partial orders used to generate filters on $\lambda$. 
\item A bar above a letter denotes that the object considered is a sequence;
usually $\bar{X}$ will be $\langle X_i:i<\zeta\rangle$, where $\zeta$
is the length $\lh(\bar{X})$ of $\bar{X}$. Sometimes our sequences will be
indexed by a set of ordinals, say $S\subseteq\lambda$, and then $\bar{X}$
will typically be $\langle X_\delta:\delta\in S\rangle$. 
\end{enumerate}

\begin{definition}
A dominating family in ${}^\lambda\lambda$ is a family $\cF\subseteq
  {}^\lambda\lambda$ such that 
\[(\forall g\in {}^\lambda\lambda)(\exists f\in \cF)(\exists\alpha<\lambda)
  (\forall\beta>\alpha)(g(\beta)<f(\beta)).\] 
The $\lambda$--dominating number $\gdl$ is defined as 
\[\gdl=\min\big\{\|\cF\|:\cF\subseteq {}^\lambda\lambda\mbox{ is a
dominating family in }{}^\lambda\lambda\;\big\}.\]
A club--dominating family in ${}^\lambda\lambda$ is a family $\cF\subseteq 
  {}^\lambda\lambda$ such that 
\[\big(\forall g\in {}^\lambda\lambda\big)\big(\exists f\in \cF\big)\big(
\{\beta<\lambda:g(\beta)\geq f(\beta)\}\mbox{ is non-stationary in }\lambda
\big).\] 
The $\cl(\lambda)$--dominating number $\gcl$ is defined as 
\[\gcl=\min\big\{\|\cF\|:\cF\subseteq {}^\lambda\lambda\mbox{ is a
$\cl(\lambda)$--dominating family in }{}^\lambda\lambda\;\big\}.\]
On $\gdl,\gcl$ see, e.g., in Cummings and Shelah \cite{CuSh:541}.
\end{definition}
\bigskip

\noindent{\bf Acknowledgment:}\quad I thank Tomek Bartoszy\'nski and Andrzej
Ros{\l}anowski for stimulating discussions. 

\section{Weakly reasonable ultrafilters}
In Definition \ref{1.5}(1) we formulate the main property of ultrafilters on
$\lambda$ which is of interest to us: being a weakly reasonable
ultrafilter. In the spectrum of all ultrafilters, weakly reasonable
ultrafilters are at the opposite end to the one occupied by normal
ultrafilters. We show that there exist (in ZFC) weakly reasonable
ultrafilters (see \ref{1.8K}) and we also give some properties of such
ultrafilters. 

\begin{definition}
\label{0.1A}
For a cardinal $\lambda$, 
\begin{enumerate}
\item[(a)]  $\ulf(\lambda)$ is the set of all ultrafilters on $\lambda$, 
\item[(b)]  $\uuf(\lambda)$ is the family of all uniform ultrafilters on
  $\lambda$, 
\item[(c)] if $D$ is a filter on $\lambda$ and $f\in {}^\lambda\mu$, then 
\[D/f\stackrel{\rm def}{=}\{A\subseteq\mu:f^{-1}(A)\in D\}\]
(usually $\mu = \lambda$). 
\end{enumerate}
\end{definition}

\begin{definition}
\label{0.1C}
Assume $D$ is an ultrafilter on $\lambda$.
\begin{enumerate}
\item If $E$ is an equivalence relation on $\lambda$, then $f_E\in
  {}^\lambda\lambda$ is defined by 
\[f_E(\alpha)=\otp\Big(\{\beta<\alpha:\beta=\min(\beta/E)<\min(\alpha/E)
\}\Big),\]
and $D/E$ is $D/f_E$.
\item  For a club $C$ of $\lambda$ let $E_C$ be the following equivalence
  relation on $\lambda$: 
\[\alpha E_C \beta\quad\mbox{ iff }\quad (\forall\gamma\in C)(\alpha<\gamma
  \ \Leftrightarrow\ \beta <\gamma),\]
and let $D/C$ be $D/E_C$. 
\item $\cFl$ is the family of all non-decreasing unbounded functions
  from $\lambda$ to $\lambda$. 
\end{enumerate}
\end{definition}

\begin{observation}  
\label{0.1aB}
Assume that $\lambda$ is a regular cardinal, $D\in \ulf(\lambda)$. 
\begin{enumerate}
\item If $f:\lambda\longrightarrow\lambda$, then $D/f\in\ulf(\lambda)$.  
\item If $f\in\cFl$ and $D$ is uniform, then also $D/f$ is a uniform
  ultrafilter on $\lambda$.   
\item If $C$ is a club of $\lambda$ and $\langle\delta_\xi:\xi<\lambda
  \rangle$  is the increasing enumeration of $C$, then for a set $A\subseteq 
  \lambda$, 
\[A\in D/C\quad\mbox{ if and only if }\quad \bigcup\big\{[\delta_\xi,
  \delta_{\xi+1}): \xi\in A\big\}\in D.\] 
\end{enumerate}
\end{observation}

\begin{definition}
\label{1.5}
Let $D$ be a uniform ultrafilter on $\lambda$. 
\begin{enumerate}
\item We say that $D$ is {\em weakly reasonable\/} if for every $f\in\cFl$
there is a club $C$ of $\lambda$ such that  
\[\bigcup\{[\delta,\delta + f(\delta)):\delta \in C\}\notin D.\]
\item We define a game $\Game_D$ between two players, Odd and Even, as
follows. A play of $\Game_D$ lasts $\lambda$ steps and during a play an
increasing continuous sequence $\bar{\alpha}=\langle\alpha_i:i<\lambda
\rangle\subseteq\lambda$ is constructed. The terms of $\bar{\alpha}$ are
chosen successively by the two players so that Even chooses the $\alpha_i$
for even $i$ (including limit stages $i$ where she has no free choice) and
Odd chooses $\alpha_i$ for odd $i$.

Even wins the play if and only if  
\[\bigcup\{[\alpha_{2i+1},\alpha_{2i+2}):i<\lambda\}\in D.\]
\end{enumerate}
\end{definition}

\begin{observation}
\label{easyob}
Let $D\in\uuf(\lambda)$. Then the following conditions are equivalent:
\begin{enumerate}
\item[(A)] $D$ is weakly reasonable,
\item[(B)] for every increasing continuous sequence $\langle\delta_\xi:\xi< 
\lambda\rangle\subseteq\lambda$ there is a club $C^*$ of $\lambda$ such that 
\[\bigcup\big\{[\delta_\xi,\delta_{\xi+1}):\xi\in C^*\big\}\notin D,\] 
\item[(C)] for every club $C$ of $\lambda$ the quotient $D/C$ does not
  extend the filter generated by clubs of $\lambda$.
\end{enumerate}
\end{observation}

\begin{proposition}
\label{1.7} 
Assume $D\in\uuf(\lambda)$. 
\begin{enumerate}
\item If $\lambda$ is strongly inaccessible and Odd has a winning strategy
  in $\Game_D$, then $D$ is not weakly reasonable.
\item If $D$ is not weakly reasonable, then Odd has a winning strategy in
  the game $\Game_D$. 
\item In part (1) instead ``$\lambda$ is strongly inaccessible", it suffices 
to assume $\diamondsuit^*_\lambda$.
\end{enumerate}
\end{proposition}

\begin{proof}
(1)\qquad Suppose towards contradiction that $\lambda$ is strongly
inaccessible, Odd has a winning strategy $\st$ in the game $\Game_D$ but
$D$ is weakly reasonable. By induction on $\varepsilon<\lambda$ choose an
increasing continuous sequence $\langle N_\varepsilon:\varepsilon<\lambda
\rangle$ of elementary submodels of $\cHl$ so that for each $\varepsilon$: 
\begin{enumerate}
\item[(a)] $N_\varepsilon\prec(\cHl,\in,<^*)$, $\|N_\varepsilon\| <
  \lambda$, $N_\varepsilon \cap \lambda \in \lambda$,
\item[(b)] ${}^\vare N_{\vare+1}\subseteq N_{\vare+1}$, 
\item[(c)] $\langle N_\zeta:\zeta\le\varepsilon\rangle\in
  N_{\varepsilon+1}$,  
\item[(d)] $\st,\lambda,D$ belong to $N_0$. 
\end{enumerate}
Let $\delta_\varepsilon=N_\varepsilon\cap\lambda$ (for $\varepsilon<
\lambda$). Thus $\langle\delta_\varepsilon:\varepsilon<\lambda\rangle$
is an increasing continuous sequence of limit ordinals. Let
$f(\alpha)=\delta_{\alpha+1}$ for $\alpha<\lambda$, so $f\in\cFl$. 

Since $D$ is a weakly reasonable ultrafilter, there is a club $C$ of
$\lambda$ such that 
\[\bigcup\{[\delta,\delta+f(\delta)):\delta\in C\}\notin D.\]
Let 
\[C^*=\big\{\varepsilon\in C:\varepsilon=\delta_\varepsilon\mbox{ is a limit  
  ordinal }\big\}\]
(it is a club of $\lambda$). Then for $\varepsilon\in C^*$ we have 
$[\delta_\vare,\delta_{\vare+1})\subseteq [\vare,\vare+f(\vare))$ and hence 
\[\bigcup\{[\delta_\varepsilon,\delta_{\vare+1}):\vare\in C^*\}\notin D.\]
Let us define a strategy $\st'$ for Even in the game $\Game_D$ as
follows. For an even ordinal $i<\lambda$, in the $i$-th move of a play, if
$\langle\alpha_j:j<i\rangle$ has been played so far then Even plays
\[\alpha_i=\left\{\begin{array}{ll}
\sup\{\alpha_j:j<i\}&\mbox{ if $i$ is limit,}\\
\min\{\vare\in C^*:(\forall j<i)(\alpha_j<\vare)\}&\mbox{ otherwise.}
\end{array}\right.\] 
Now consider a play $\langle\alpha_i:i<\lambda\rangle$ in which Even uses
the strategy $\st'$ and Odd plays according to $\st$. Then for each
$i<\lambda$ we have $\alpha_{2i}\in C^*$ and thus $\alpha_{2i}=\delta_{
\alpha_{2i}}\in N_{\alpha_{2i}+1}$, and also $\{\alpha_j:j<2i\}\subseteq
\alpha_{2i}\subseteq N_{\alpha_{2i}+1}$. Since the model $N_{\alpha_{2i}+1}$
is closed under forming sequences of length $\alpha_{2i}+1$ (by (b)), we
conclude that $\langle\alpha_j:j\leq 2i\rangle \in N_{\alpha_{2i}+1}$. Since
$\st\in N_0\prec N_{\alpha_{2i}+1}$, clearly $\alpha_{2i+1}\in 
N_{\alpha_{2i}+1}\cap\lambda$ and therefore $\alpha_{2i+1}<
\delta_{\alpha_{2i}+1}$. Hence
\[\bigcup\{[\alpha_{2i},\alpha_{2i+1}):i<\lambda\}\subseteq\bigcup
\{[\delta_{\alpha_{2i}},\delta_{\alpha_{2i}+1}):i<\lambda\}\subseteq
\bigcup\{[\delta_\vare,\delta_{\vare +1}):\vare\in C^*\}\notin D.\] 
But $\st$ is a winning strategy for Odd, so he wins the play and
$\bigcup\{[\alpha_{2i+1},\alpha_{2i+2}):i<\lambda\}\notin D$, a
contradiction. 
\medskip

\noindent (2)\quad Suppose that $D\in\uuf(\lambda)$ is not weakly
reasonable. Then we may find $f\in\cFl$ such that for every club $C$ of 
$\lambda$ we have 
\[\bigcup\{[\delta,\delta+f(\delta)):\delta\in C\}\in D.\]
Let $\st$ be a strategy of Odd in $\Game_D$ which instructs him to play as
follows. For an odd ordinal $i=i_0+1<\lambda$, in the $i$-th move of a play,
if $\langle\alpha_j:j\leq i_0\rangle$ has been played so far, then Odd plays 
$\alpha_i=\alpha_{i_0}+f(\alpha_{i_0})+1$. 

We claim that $\st$ is a winning strategy for Odd (in $\Game_D$). To this
end suppose that $\langle\alpha_j:j<\lambda\rangle\subseteq\lambda$ is a
result of a play of $\Game_D$ in which Odd uses the strategy $\st$. Let 
$C'=\{\alpha_i:i<\lambda\mbox{ is limit }\}$ -- it is a club of $\lambda$,
so by the choice of $f$ we have  
\[\bigcup\{[\delta,\delta+f(\delta)):\delta\in C'\}\in D.\]
Since $\bigcup\{[\delta,\delta+f(\delta)):\delta\in C'\}\subseteq
\bigcup\{[\alpha_{2i},\alpha_{\alpha_{2i+1}}): i<\lambda\}$ we may now
  conclude that Odd indeed wins the play.  
\end{proof}

\begin{lemma}
\label{1.7Q}
Suppose that $\lambda$ is a regular uncountable cardinal, $D\in\uuf(
\lambda)$ is a weakly reasonable ultrafilter and $\langle\beta_i:i<\lambda
\rangle$ is an increasing continuous sequence of ordinals below $\lambda$.
Then there is an increasing continuous sequence $\langle\delta_\xi:\xi<
\lambda\rangle\subseteq\lambda$ consisting of limit ordinals and such that  
\[\bigcup\{[\beta_{\delta_{2\xi+1}},\beta_{\delta_{2\xi+2}}):\xi<\lambda\}
  \in D.\]   
\end{lemma}

\begin{proof}
It follows from \ref{easyob} that we may find a club $C^*$ of $\lambda$ such 
that all members of $C^*$ are limit ordinals and $\bigcup\big\{[\beta_\xi,
\beta_{\xi+1}):\xi\in C^*\big\}\notin D$. Let $C^+=C^*\cup\{\xi+1:\xi\in
C^*\}$ (clearly it is a club of $\lambda$) and let $\langle\delta_\xi:\xi<
\lambda\rangle$ be the increasing enumeration of $C^+$. Note that
$C^*=\{\delta_\xi:\xi<\lambda\mbox{ is even }\}$ and, for an even ordinal
$\xi<\lambda$, $\delta_{\xi+1}=\delta_\xi+1$. Hence
\[\begin{array}{r}
\bigcup\big\{[\beta_{\delta_\xi},\beta_{\delta_{\xi+1}}):\xi<\lambda
\mbox{ is even }\big\}=\bigcup\big\{[\beta_{\delta_\xi},
\beta_{\delta_\xi+1}):\xi<\lambda\mbox{ is even }\big\}=\qquad\\
\bigcup\big\{[\beta_\zeta,\beta_{\zeta+1}):\zeta\in C^*\big\}\notin D
  \end{array}\] 
Consequently, $\bigcup\big\{[\beta_{\delta_\xi},\beta_{\delta_{\xi+1}}):\xi
<\lambda\mbox{ is odd }\big\}\in D$. 
\end{proof}

\begin{theorem}
\label{1.8}
If $\lambda$ is a regular uncountable cardinal and $D\in\uuf(\lambda)$ is
weakly reasonable, then $D$ is a regular ultrafilter. 
\end{theorem}

\begin{proof}
Using Lemma \ref{1.7Q} we may choose by induction on $\vare<\lambda$ a
sequence $\langle\bar{\delta}^\vare:\vare<\lambda\rangle$ so that 
\begin{enumerate}
\item[(a)] $\bar{\delta}^\varepsilon=\langle\delta^\varepsilon_i:i<
\lambda\rangle$ is an increasing continuous sequence of non-successor
ordinals below $\lambda$, $\delta^\varepsilon_0 = 0$, 
\item[(b)] the set $A_\varepsilon\stackrel{\rm def}{=}\bigcup\big\{[
\delta^\varepsilon_{2i+1},\delta^\varepsilon_{2i+2}):i<\lambda\}$ belongs to
$D$,  
\item[(c)] if $\zeta<\varepsilon$, $i<\lambda$, then $\delta^\varepsilon_i
  \in \{\delta^\zeta_j:j<\lambda$ is a limit ordinal or zero $\}$. 
\end{enumerate}
For $\vare<\lambda$ let $f_\vare:A_\vare\longrightarrow\lambda$ be such that 
\[\alpha\in [\delta^\varepsilon_{2i+1},\delta^\varepsilon_{2i+2})\quad
\Rightarrow\quad f_\varepsilon(\alpha)=\delta^\varepsilon_{2i+1}.\]
Note that
\begin{enumerate}
\item[$(\otimes)$] if $\zeta<\vare<\lambda$, $\alpha\in A_\zeta\cap
  A_\vare$, then $f_\vare(\alpha)<f_\zeta(\alpha)$. 
\end{enumerate}
[Why? Let $f_\zeta(\alpha)=\delta^\zeta_{2i+1}$ (so $\alpha\in
[\delta^\zeta_{2i+1},\delta^\zeta_{2i+2})$). It follows from (c) that
$f_\vare(\alpha)\in \{\delta^\zeta_j:j<\lambda$ is a limit ordinal or zero
$\}$ and hence (since also $f_\vare(\alpha)\leq \alpha$) we may conclude 
that $f_\vare(\alpha)<f_\zeta(\alpha)$. ]

For $\alpha<\lambda$, let $w_\alpha=\{\vare<\lambda:\alpha\in A_\vare\}$. It 
follows from $(\otimes)$ that (for every $\alpha<\lambda$) the sequence
$\langle f_\vare(\alpha):\vare\in w_\alpha\rangle$ is strictly decreasing,
so necessarily each $w_\alpha$ is finite. Since $A_\vare\in D$ for each
$\vare<\lambda$ (by (b)), we have shown the regularity of $D$.
\end{proof}

\begin{theorem}
\label{1.8K}
Let $\lambda>\aleph_0$ be a regular cardinal. Then there is a uniform weakly
reasonable ultrafilter $D$ on $\lambda$. 
\end{theorem}

\begin{proof}
Let $\{f_\vare:\vare<\gdl\}\subseteq {}^\lambda\lambda$ be a dominating 
family and for $\vare<\gdl$ let $C_\vare$ be a club of $\lambda$ such that
members of $C_\vare$ are limit ordinals and 
\[(\forall\delta\in C_\vare)(\forall\alpha<\delta)(f_\vare(\alpha)<
\delta).\]
Let $\langle\alpha_{\vare,i}:i<\lambda\rangle$ be the increasing enumeration
of $C_\vare$.  

By induction on $\vare$ we will choose sets $E_\vare,A_\vare$ so that for
each $\vare<\gdl$: 
\begin{enumerate}
\item[(a)]  $A_\vare$ is an unbounded subset of $\lambda$ and $E_\vare
\subseteq C_\vare$ is a club of $\lambda$,  
\item[(b)]  $A_\vare\cap\bigcup\big\{[\alpha_{\vare,\gamma},\alpha_{\vare,
    \gamma+1}): \gamma\in E_\vare\}=\emptyset$, 
\item[(c)]  if $n<\omega$, $\zeta_0<\ldots<\zeta_{n-1}<\vare$, then
  $\|A_\vare\cap \bigcap\limits_{i<n}A_{\zeta_i}\|=\lambda$.
\end{enumerate}
So suppose that we have chosen $A_\zeta,E_\zeta$ for $\zeta<\vare<\gdl$
so that the respective reformulations of (a)--(c) hold true. For a finite
sequence $\bar{\zeta}=\langle\zeta_i:i<n\rangle$ of ordinals below $\vare$
let $A^{\bar{\zeta}}=\bigcap\limits_{i<n}A_{\zeta_i}$ (note that
$\|A^{\bar{\zeta}}\|=\lambda$ by the demand in (c)). Let
$g_{\bar{\zeta}}^\vare\in {}^\lambda\lambda$ be such that   
\begin{enumerate}
\item[$(\oplus)$]  if $\alpha_{\vare,i}\leq\alpha<\alpha_{\vare,i+1}$, then
$g_{\bar{\zeta}}^\vare(\alpha)=\min\{\delta>\alpha_{\vare,i+1}:[
\alpha_{\vare,i+1},\delta)\cap A^{\bar{\zeta}}\neq\emptyset\}$. 
\end{enumerate}
The family $\{g_{\bar{\zeta}}^\vare:\bar{\zeta}\in {}^{\omega>}\vare\}$ is a   
subset of ${}^\lambda\lambda$ of cardinality $\leq |\vare|+\aleph_0<\gdl$, 
so it cannot be a dominating family. Therefore we may pick a function
$h_\vare\in {}^\lambda\lambda$ such that 
\[(\forall\bar{\zeta}\in {}^{\omega>}\vare)(\exists^\lambda\alpha<\lambda)
(g_{\bar{\zeta}}^\vare(\alpha)<h_\vare(\alpha)).\] 
Put 
\[\begin{array}{l}
E_\vare=\{\delta<\lambda:\delta=\alpha_{\vare,\delta}\mbox{ is a limit
  ordinal and }(\forall\alpha<\delta)(h_\vare(\alpha)<\delta)\}\quad
\mbox{ and}\\  
A_\vare=\bigcup\big\{[\alpha_{\vare,\gamma+1},\alpha_{\vare,\delta}):\gamma
  <\delta\mbox{ are successive members of }E_\varepsilon\}.
  \end{array}\]
It should be clear that $E_\vare,A_\vare$ satisfy demands (a), (b). 

Let us argue that also condition (c) holds true. Let $\bar{\zeta}\in
{}^{\omega>} \vare$ and we shall prove that $A_\vare\cap A^{\bar{\zeta}}$ is 
unbounded in $\lambda$. By the choice of $h_\vare$, the set
$B=\{\alpha<\lambda:g_{\bar{\zeta}}^\vare(\alpha)<h_\vare(\alpha)\}$ is of
cardinality $\lambda$.  Let us fix for a moment $\alpha \in B$ and let
$i<\lambda$ be such that $\alpha_{\vare,i}\leq \alpha<\alpha_{\vare,i+1}$.
Let $\sup(E_\vare\cap\alpha_{\varepsilon,i+1})=\gamma=\alpha_{\vare,\gamma}$
and $\min(E_\vare\setminus\alpha_{\vare,\gamma+1})=\delta=\alpha_{\vare,
\delta}$. Then $\gamma,\delta$ are successive members of $E_\vare$ and 
\[\gamma\leq\alpha_{\vare,i}\leq\alpha<\alpha_{\vare,i+1}<\delta.\]
Hence (by the definition of $E_\vare$ and by $\alpha\in B$) we get 
\[[\alpha_{\varepsilon,i+1},g_{\bar{\zeta}}^\vare(\alpha))\subseteq
[\alpha_{\vare,i+1},h_\vare(\alpha))\subseteq [\alpha_{\vare,\gamma+1},
\alpha_{\vare,\delta})\subseteq A_\varepsilon.\]
It follows from $(\oplus)$ that $[\alpha_{\vare,i+1},g_{\bar{\zeta}}^\vare
(\alpha))\cap A^{\bar{\zeta}}\neq\emptyset$, and consequently $A_\vare\cap
A^{\bar{\zeta}}\setminus \alpha\neq\emptyset$. Since $\|B\|=\lambda$ we may
now easily conclude that $\|A_\vare\cap A^{\bar{\zeta}}\|=\lambda$, showing
that $A_\vare,E_\vare$ are as required. 
\medskip

After the construction is carried out (and we have the sequence $\langle
E_\vare,A_\vare:\vare<\gdl\rangle$) we may find a uniform ultrafilter $D$ on
$\lambda$ such that $\{A_\vare:\vare<\gdl\}\subseteq D$ (remember the demand
in (c)). We claim that $D$ is weakly reasonable. To this end suppose that
$C$ is a club of $\lambda$ and $\langle\delta_\xi:\xi<\lambda\rangle
\subseteq\lambda$ is the increasing enumeration of $C$. By the choice of
$f_\vare,C_\vare$ (for $\vare<\gdl$) we may find $\vare<\gdl$ and
$j_0<\lambda$ such that  
\[(\forall i\geq j_0)(\|\;[\alpha_{\vare,i},\alpha_{\vare,i+1})\cap C\;\|
  >2).\] 
Let 
\[C^*=\{\gamma\in E_\vare\cap C\setminus j_0:\gamma=\alpha_{\vare,\gamma}=
\delta_\gamma\mbox{ is a limit ordinal }\}\]
(it is a club of $\lambda$). Since for $\gamma\in C^*$ we have that
$\alpha_{\vare,\gamma}=\delta_\gamma<\delta_{\gamma+1}<\alpha_{\vare,
\gamma+1}$ we may easily conclude from (b) that 
\[\bigcup\big\{[\delta_\gamma,\delta_{\gamma+1}):\gamma\in C^*\big\}\notin
  D,\] 
completing the proof (remember \ref{easyob}).
\end{proof}

\section{More reasonable ultrafilters}
In this section we propose a property of ultrafilters stronger than being 
weakly reasonable (see Definition \ref{1.3}(5)). We believe that the notion
of {\em very reasonable ultrafilters\/} is the right re-interpretation of
being a $P$--point in the setting of ``very non-normal ultrafilters'' on an 
uncountable regular cardinal $\lambda$. We start with describing a forcing
notion $\bql$ which motivated our choice of {\em generating systems\/} of 
\ref{1.3}.  

Like before, $\lambda$ is always assumed to be an uncountable regular
cardinal.  

\begin{definition}
\label{1.1}
We define a forcing notion $\bql$ as follows.\\
{\bf A condition in $\bql$} is a tuple $p=(\gamma^p,C^p,\langle Z^p_\delta: 
\delta\in C^p\rangle,\langle d^p_\delta:\delta \in C^p \rangle)$ such that 
\begin{enumerate}
\item[(i)]   $\gamma^p<\lambda$, $C^p$ a club of $\lambda$ consisting of
  limit ordinals only, and for $\delta\in C^p$:
\item[(ii)]  $Z^p_\delta=\Big[\delta,\min\big(C^p\setminus(\delta+1)\big)
\Big)$ and 
\item[(iii)] $d^p_\delta\subseteq\cP(Z^p_\delta)$ is a proper ultrafilter on 
  $Z^p_\delta$. 
\end{enumerate}
{\bf The order $\leq_{\bql}=\leq$ of $\bql$} is given by \\
$p\leq_{\bql} q$\qquad if and only if \\
\begin{enumerate}
\item[(a)] $\gamma^p\leq\gamma^q$, $C^p\cap\gamma^p\subseteq C^q\subseteq
  C^p$ and 
\item[(b)] if $\delta<\vare$ are successive members of $C^q$ (so
  $Z^q_\delta=[\delta,\vare)$), then 
\[\big(\forall A\in d^q_\delta\big)\big(\exists\zeta\in C^p\cap
  [\delta,\vare)\big)\big(A\cap Z^p_\zeta\in d^p_\zeta\big).\]
\end{enumerate}
\end{definition}

\begin{remark} The forcing notion $\bql$ can be represented according to
  the framework of \cite[\S B.5]{RoSh:777}.
\end{remark}

\begin{proposition}
\label{1.2} 
\begin{enumerate}
\item $\bql$ is a partial order, $\|\bql\|=2^{2^{<\lambda}}$.
\item If $p,q\in\bql$, $p\leq q$, $\delta<\vare$ are two successive members 
of $C^p$, $\delta,\vare\in C^q$, then $Z^q_\delta=Z^p_\delta$ and
$d^q_\delta=d^p_\delta$.
\item  $\bql$ is $({<}\lambda)$---complete (so it does not add bounded
  subsets of $\lambda$). 
\item  If $p\in\bql$, $A\subseteq\lambda$, then there is a condition $q\in
  \bql$ stronger than $p$ and such that 
\[\mbox{either}\quad (\forall\delta\in C^q\setminus\gamma^p)(A\cap
  Z^q_\delta\in d^q_\delta)\quad\mbox{ or }\quad (\forall\delta\in
  C^q\setminus \gamma^p)(A\cap Z^q_\delta\notin d^q_\delta).\]  
\end{enumerate}
\end{proposition}

\begin{proof}
(1), (2)\quad Straightforward.
\medskip

\noindent (3)\quad Assume that $\delta<\lambda$ is a limit ordinal and a
sequence $\langle p_i:i<\delta \rangle\subseteq\bql$ is
$\le_{\bql}$---increasing. Let $E$ be a uniform ultrafilter on $\delta$.
Let us put:   
\begin{itemize}
\item $\gamma=\sup\{\gamma^{p_i}:i<\delta\}$, $C=\bigcap\limits_{i<\delta}
C^{p_i}$, and for $\alpha \in C$ let
\item $Z_\alpha=\big[\alpha,\min(C\setminus(\alpha+1))\big)$ and 
\item $d_\alpha=\big\{A\subseteq Z_\alpha:\{i<\delta: A\cap
  Z^{p_i}_\alpha\in d^{p_i}_\alpha\} \in E\big\}$.
\end{itemize}
It is easy to check that $p=(\gamma,C,\langle Z_\alpha:\alpha\in C\rangle,
\langle d_\alpha:\alpha\in C\rangle)$ belongs to $\bql$ and that it is a
condition stronger than all $p_i$ (for $i<\delta$).
\medskip

\noindent (4) Let $p\in\bql$ and $A\subseteq\lambda$. Just for simplicity we
may assume that $\gamma^p\in C^p$ (as we may always increase $\gamma^p$). Put 
\[Y\stackrel{\rm def}{=}\{\alpha\in C^p:A \cap Z^p_\alpha\in
d^p_\alpha\}\]
and let us consider two cases.

{\sc Case 1:}\quad $Y$ is unbounded in $\lambda$.\\
Then we may choose an increasing continuous sequence $\langle\delta_i:
i<\lambda\rangle\subseteq C^p$ such that $\delta_0=\gamma^p$ and
$\big(\forall i<\lambda\big)\big([\delta_i,\delta_{i+1})\cap Y\ne\emptyset
\big)$. Put
\begin{itemize}
\item $\gamma=\gamma^p$, $C=\{\delta_i:i<\lambda\}\cup(C^p\cap\gamma^p)$,  
\item if $\alpha\in C^p\cap\gamma^p$, then $Z_\alpha=Z^p_\alpha$ and
  $d_\alpha=d^p_\alpha$, 
\item if $\alpha=\delta_i$, $i<\lambda$, then $Z_\alpha=[\delta_i,
\delta_{i+1})$ and 
\[d_\alpha=\big\{B\subseteq Z_\alpha:B\cap Z^p_{\min(Y\setminus\alpha)}\in 
d^p_{\min(Y\setminus\alpha)}\big\}.\]  
\end{itemize}
It is straightforward to verify that $q=(\gamma,C,\langle Z_\alpha:\alpha\in
C\rangle, \langle d_\alpha:\alpha\in C\rangle)\in\bql$ is a
condition stronger than $p$ and it is also clear that $(\forall\alpha\in
C\setminus\gamma^p)(A\cap Z_\alpha\in d_\alpha)$.  

{\sc Case 2:}\quad $Y$ is bounded in $\lambda$.\\
Then the set $\lambda\setminus Y$ is unbounded, so we may apply the
construction of $q$ from Case 1 replacing $Y$ by its complement
$\lambda\setminus Y$. It should be clear that the condition $q$ which we
will get then satisfies $(\forall\alpha\in C^q\setminus\gamma^p)(A\cap
Z^q_\alpha\notin d_\alpha^q)$.   
\end{proof}

\begin{remark}
\label{userem}
The following discussion presents our motivations for the definitions and
concepts presented later in this section.

Suppose that $G\subseteq\bql$ is a generic filter over $\bV$. In $\bV[G]$ we
define $C=\bigcup\{C^p\cap\gamma^p:p\in G\}$ and for $\alpha\in C$ we let
$d_\alpha=d^p_\alpha$ for some (equivalently: all) $p\in G$ such that
$\alpha<\gamma^p$ and $C^p\cap (\alpha,\gamma^p)\neq\emptyset$. Then $C$ is
a club of $\lambda$ and (for $\alpha\in C$) $d_\alpha$ is an ultrafilter on
$[\alpha,\min(C\setminus(\alpha+1)))$. Let 
\[D=\big\{A\in\cP(\lambda)^\bV: \big(\exists\vare<\lambda\big)\big(\forall
\alpha >\vare\big)\big(A\cap [\alpha,\min(C\setminus(\alpha+1)))\in d_\alpha
  \big)\big\}.\]
It follows from \ref{1.2}(4) that $D$ is an ultrafilter on the Boolean
  algebra $\cP(\lambda)^\bV$.

Let $\name{D}$ be a $\bql$--name for the $D$ defined as above. Note that if
$p\in\bql$, $A\subseteq\lambda$ and $(\exists\vare<\lambda)(\forall\delta\in
C^p\setminus\vare)(A\cap Z^p_\delta\in d^p_\delta)$, then $p\forces_{\bql}
\mbox{`` }A\in\name{D}\mbox{ ''}$. Plainly, the family $\{A\subseteq\lambda: 
p\forces_{\bql}\mbox{`` }A\in\name{D}\mbox{ ''}\}$ is a uniform filter on
$\lambda$, and, of course, for a generic filter $G\subseteq\bql$ over $\bV$,  
\[\name{D}^G=\bigcup\Big\{\{A\subseteq\lambda: p\forces_{\bql}\mbox{`` }A\in
\name{D}\mbox{ ''}\}:p\in G\Big\}.\]
\end{remark}

\begin{definition}
\label{1.3} 
\begin{enumerate}
\item We define a forcing notion $\bqz$ as follows.\\
{\bf A condition in $\bqz$} is a tuple $p=(C^p,\langle Z^p_\delta: 
\delta\in C^p\rangle,\langle d^p_\delta:\delta \in C^p \rangle)$ such that 
$(0,C^p,\langle Z^p_\delta:\delta\in C^p\rangle,\langle d^p_\delta:\delta
\in C^p \rangle)\in\bql$; \\
{\bf The order $\leq_{\bqz}=\leq$ of $\bqz$} is inherited from $\bql$ in a
natural way.
\item We define a relation $\leq^*_{\bqz}=\leq^*$ on $\bqz$ as follows: 

$p\leq^* q$\quad if and only if \quad for some $\alpha<\lambda$ we have 
\[(C^p\setminus\alpha,\langle Z^p_\delta:\delta\in C^p\setminus\alpha
\rangle,\langle d^p_\delta:\delta\in C^p\setminus\alpha\rangle)\leq_{\bqz}
(C^q\setminus\alpha,\langle Z^q_\delta:\delta\in C^q\setminus\alpha
\rangle,\langle d^q_\delta:\delta\in C^q\setminus\alpha\rangle).\]
\item For a condition $q\in\bqz$ we let 
\[\fil(q)\stackrel{\rm def}{=}\big\{A\subseteq\lambda: (\exists\vare<
\lambda)(\forall\delta\in C^q\setminus\vare)(A\cap Z^q_\delta\in
d^q_\delta)\big\},\]  
and for a set $G^*\subseteq\bqz$ we let $\fil(G^*)\stackrel{\rm def}{=}
\bigcup\{\fil(p):p\in G^*\}$. We also define a binary relation $\leq^0$ on
$\bqz$ by 
\begin{center}
$p\leq^0 q$\quad if and only if\quad $\fil(p)\subseteq\fil(q)$.
\end{center}
\item We say that an ultrafilter $D$ on $\lambda$ is {\em reasonable\/} if
it is weakly reasonable (see \ref{1.5}(1)) and there is a
directed (with respect to $\leq^0$) set $G^*\subseteq\bqz$ such that
$D=\fil(G^*)$. The family $G^*$ may be called {\em the generating system for
  $D$}.    
\item An ultrafilter $D$ on $\lambda$ is said to be {\em very reasonable} if 
it is weakly reasonable and there is a $({<}\lambda^+)$--directed (with
respect to $\leq^0$) set $G^*\subseteq\bqz$ such that $D=\fil(G^*)$. 
\end{enumerate}
\end{definition}

\begin{remark}
Note that $\|\fil(p)\|=2^\lambda$ for each $p\in\bqz$. Thus even if
$D=\fil(G^*)$ for some small generating system $G^*\subseteq\bqz$, the
minimal number of generators for $D$ as a filter may be $2^\lambda$. 
\end{remark}

\begin{observation}
\begin{enumerate}
\item If $p\leq^*_\bqz q$, then $\fil(p)\subseteq\fil(q)$ (so $p\leq^0 q$).  
\item If a set $G^*\subseteq\bqz$ is directed with respect to $\leq^0$, then
$\fil(G^*)$ is a filter of subsets of $\lambda$ containing all co-bounded
  subsets of $\lambda$. 
\end{enumerate}
\end{observation}

\begin{definition}
\label{eproduct} 
Suppose that
\begin{enumerate}
\item[(a)] $X$ is a non-empty set and $e$ is an ultrafilter on $X$,
\item[(b)] $d_x$ is an ultrafilter on a set $Z_x$ (for $x\in X$). 
\end{enumerate}
We let 
\[\bigoplus\limits^e_{x\in X} d_x=\big\{A\subseteq\bigcup\limits_{x\in X}
Z_x: \{x\in X:Z_x\cap A\in d_x\}\in e\big\}.\]
(Clearly, $\bigoplus\limits^e_{x\in X} d_x$ is an ultrafilter on
$\bigcup\limits_{x\in X} Z_x$.)
\end{definition}

\begin{proposition}
\label{tfcae}
Let $p,q\in\bqz$. Then the following are equivalent:
\begin{enumerate}
\item[(a)]  $p\leq^0 q$,
\item[(b)] there is $\vare<\lambda$ such that  
\[\big(\forall\alpha\in C^q\setminus\vare\big)\big(\forall A\in d^q_\alpha 
\big)\big(\exists \beta\in C^p\big)\big(A\cap Z^p_\beta\in d^p_\beta
\big),\]    
\item[(c)] there is $\vare<\lambda$ such that 

if $\alpha\in C^q\setminus\vare$, $\beta_0=\sup\big(C^p\cap(\alpha+1)\big)$, 
$\beta_1=\min\big( C^p\setminus\min(C^q\setminus(\alpha+1))\big)$,

then there is an ultrafilter $e$ on $[\beta_0,\beta_1)\cap C^p$ such that 
\[d^q_\alpha=\big\{A\cap Z^q_\alpha:A\in\bigoplus\limits^e\{d^p_\beta:
\beta\in [\beta_0,\beta_1)\cap C^p\}\big\}.\] 
\end{enumerate}
\end{proposition}

\begin{proof}
(a) $\Rightarrow$ (b)\quad Assume towards contradiction that $p\leq^0 q$,
but (b) fails. Then we may pick a sequence $\langle \alpha_\xi,A_\xi:\xi<
\lambda\rangle$ such that for each $\xi<\lambda$,   
\begin{enumerate}
\item[(i)]   $\alpha_\xi\in C^q$, $A_\xi\in d_\xi^q$, 
\item[(ii)]  if $\xi<\zeta<\lambda$, $\beta\in C^p\cap\min\big(C^q\setminus
(\alpha_\xi+1)\big)$, then $\min\big(C^p\setminus (\beta+1)\big)<
  \alpha_\zeta$,  
\item[(iii)] $\big(\forall\beta\in C^p\big)\big(A_\xi\cap Z^p_\beta\notin 
  d^p_\beta\big)$. 
\end{enumerate}
It follows from (ii) that for every $\beta\in C^p$ there is at most one
$\xi<\lambda$ such that $Z^p_\beta\cap Z^q_\xi\neq\emptyset$. Also if
$\beta\in C^p$ and $Z^p_\beta\cap Z^q_\xi\in d^p_\beta$, then $(Z^q_\xi
\setminus A_\xi)\cap Z^p_\beta\in d^p_\beta$.  

Put $A=\bigcup\limits_{\xi<\lambda}A_\xi$. By what we have said above,
for all $\beta\in C^p$ we have $(\lambda\setminus A)\cap Z^p_\beta\in
d^p_\beta$, and hence $\lambda\setminus A\in\fil(p)\subseteq\fil(q)$. This
contradicts (i).  
\medskip

\noindent (b) $\Rightarrow$ (c)\quad Assume that (b) holds true as
witnessed by $\vare<\lambda$. Let $\alpha\in C^q\setminus\vare$,
$\alpha'=\min\big(C^q\setminus (\alpha+1)\big)$, $\beta_0=\sup\big(C^p\cap 
(\alpha+1)\big)$ and $\beta_1=\min(C^p\setminus\alpha')$. For $A\in
d^q_\alpha$ put   
\[w(A)=\big\{\beta\in [\beta_0,\beta_1)\cap C^p: A\cap Z^p_\beta\in
  d^p_\beta\big\}.\]
It follows from (b) that $w(A)\neq\emptyset$. Plainly $w(A\cap A')=w(A)\cap
w(A')$ for $A,A'\in d^q_\alpha$, so we may pick an ultrafilter $e$ on
$[\beta_0,\beta_1)\cap C^p$ such that $\{w(A):A\in d^q_\alpha\}\subseteq
  e$. Now it should be clear that 
\[d^q_\alpha\subseteq \big\{B\cap Z^q_\alpha:B\in\bigoplus\limits^e\{
d^p_\beta:\beta\in [\beta_0,\beta_1)\cap C^p\}\big\}\]
and (since clearly $Z^q_\alpha\in\bigoplus\limits^e\{d^p_\beta:\beta\in
[\beta_0,\beta_1)\cap C^p\}$) the set on the right-hand side is a proper
filter on $Z^q_\alpha$. Consequently the two sets are equal. 
\medskip

\noindent (c) $\Rightarrow$ (a)\quad Assume that (c) holds true as
witnessed by $\vare<\lambda$, and suppose that $A\in\fil(p)$. Pick $\vare'
<\lambda$ such that $\vare<\vare'$ and 
\[\big(\forall\beta\in C^p\setminus\vare'\big)\big(A\cap Z^p_\beta\in
d^p_\beta\big).\] 
Suppose $\alpha\in C^q\setminus\big(\min(C^p\setminus\vare')+1\big)$ and let
$\beta_0=\sup\big(C^p\cap(\alpha+1\big)$, $\beta_1=\min\big(C^p\setminus
\min(C^q\setminus(\alpha+1))\big)$. Let $e$ be an ultrafilter on
$[\beta_0,\beta_1)\cap C^p$ such that 
\[d^q_\alpha=\big\{B\cap Z^q_\alpha:B\in\bigoplus\limits^e\{d^p_\beta:\beta
\in [\beta_0,\beta_1)\cap C^p\}\big\}.\] 
Note that $\beta_0\geq\vare'$ and hence $A\cap Z^p_\beta\in d^p_\beta$ for
all $\beta\in [\beta_0,\beta_1)\cap C^p$. Consequently 
\[A\cap [\beta_0,\beta_1)\in\bigoplus\limits^e\{d^p_\beta:\beta\in
  [\beta_0,\beta_1)\cap C^p\}\] 
and therefore also 
\[A\cap Z^q_\alpha=\big(A\cap [\beta_0,\beta_1)\big)\cap Z^q_\alpha
    \in d^q_\alpha.\] 
Now we easily conclude that $A\in \fil(q)$. 
\end{proof}

\begin{definition}
Let $p\in\bqz$. Suppose that $X\in [C^p]^{\textstyle\lambda}$ and
$C\subseteq C^p$ is a club of $\lambda$ such that
\begin{quotation}
if $\alpha<\beta$ are successive elements of $C$,\\
then $|[\alpha,\beta)\cap X|=1$.
\end{quotation}
(In this situation we say that {\em $p$ is restrictable to $\langle X,C
\rangle$}.) We define {\em the restriction of $p$ to $\langle X,C
\rangle$} as an element $q=p\rest\langle X,C\rangle\in \bqz$ such that 
$C^q=C$, and if $\alpha<\beta$ are successive elements of $C$, $x\in
[\alpha,\beta)\cap X$, then $Z^q_\alpha=[\alpha,\beta)$ and $d^q_\alpha=
\{A\subseteq Z^q_\alpha: A\cap Z^p_x\in d^p_x\}$. 
\end{definition}
 
\begin{proposition}
\label{2.10}
\begin{enumerate}
\item Assume that $G^*\subseteq\bqz$ is $\leq^0$--directed and
  $\leq^0$--downward closed, $p\in G^*$, $X\in [C^p]^{\textstyle\lambda}$
  and $C\subseteq C^p$ is a club of $\lambda$ such that $p$ is restrictable
  to $\langle X,C\rangle$. If $\bigcup\limits_{x\in X}Z_x^p\in\fil(G^*)$,
  then $p\rest\langle X,C\rangle\in G^*$. 
\item  If $G^*\subseteq\bqz$ is $\leq^0$--directed and $\|G^*\|\leq\lambda$,
  then $G^*$ has a $\leq^0$--upper bound. (Hence, in particular, $\fil(G^*)$
  is not an ultrafilter.) 
\end{enumerate}
\end{proposition}

\begin{proof}
\noindent (1)\quad Suppose that $G^*,p,X,C$ are as in the assumptions and
$\bigcup\limits_{x\in X}Z_x^p\in\fil(G^*)$. Since $G^*$ is
$\leq^0$--directed (and $p\in G^*$) we may pick $r\in G^*$ such that
$p\leq^0 r$ and $\bigcup\limits_{x\in X}Z_x^p\in\fil(r)$. We are going to
show that $q\stackrel{\rm def}{=} p\rest \langle X,C\rangle \leq^0 r$ (which
will imply that $q\in G^*$ as $G^*$ is downward closed). 

Since $\bigcup\limits_{x\in X}Z_x^p\in\fil(r)$, there is $\vare<\lambda$
such that 
\[\big(\forall\alpha\in C^r\setminus\vare\big)\big(\bigcup\limits_{x\in X}
Z_x^p\cap Z^r_\alpha\in d^r_\alpha\big)\quad\mbox{and}\quad 
\big(\alpha\in C^r\setminus\vare\big)\big(\forall A\in d^r_\alpha\big)
\big(\exists\beta\in C^p\big)\big(A\cap Z^p_\beta\in d^p_\beta\big)\]
(remember \ref{tfcae}(b)). Now suppose that $\alpha\in C^r\setminus \vare$
and $A\in d^r_\alpha$. Then $\bigcup\limits_{x\in X}Z_x^p\cap A\in
d^r_\alpha$ so there is $\beta\in C^p$ such that $\bigcup\limits_{x\in X}
Z_x^p\cap A\cap Z^p_\beta\in d^p_\beta$. In particular,
$\bigcup\limits_{x\in X}Z_x^p\cap Z^p_\beta\cap A\neq\emptyset$, so
necessarily $\beta\in X$. Let $\beta_0<\beta_1$ be the successive elements 
of $C$ such that $\beta_0\leq\beta <\beta_1$. Easily 
\[Z^p_\beta\cap A=\bigcup\limits_{x\in X}Z_x^p\cap Z^p_\beta\cap A\in
d^p_\beta\] 
and thus $A\cap Z^q_{\beta_0}\in d^q_{\beta_0}$. Thus we have shown that 
\begin{quotation}
if $\alpha\in C^r\setminus\vare$ and $A\in d^r_\alpha$,\\
then there is $\beta_0\in C^q$ such that $A\cap Z^q_{\beta_0}\in
d^q_{\beta_0}$. 
\end{quotation}
Consequently, $q\leq^0 r$ (remember \ref{tfcae}). 
\medskip

\noindent (2)\quad Let $\langle p_\xi:\xi<\lambda\rangle$ list (with
possible repetitions) all members of $G^*$. For $\xi<\lambda$ let 
$C_\xi=\big\{\delta<\lambda:\delta=\sup(\delta\cap C^{p_\xi})\big\}$ 
(it is a club of $\lambda$), and for $\xi,\zeta<\lambda$ let $\vare(\{
\xi,\zeta\})<\lambda$ be such that  
\smallskip

\noindent if $p_\xi\leq^0 p_\zeta$, then 
\[\big(\forall \alpha\in C^{p_\zeta}\setminus\vare(\{\vare,\zeta\})\big)
\big(\forall A\in d^{p_\zeta}_\alpha\big)\big(\exists \beta\in C^{p_\xi}
\big)\big(A\cap Z^{p_\xi}_\beta\in d^{p_\xi}_\beta)\big)\]

\noindent (remember \ref{tfcae}). Let 
\[C^*=\big\{\delta<\lambda:\delta\mbox{ is limit and }\{p_\xi:\xi<\delta\}
\mbox{ is $\leq^0$--directed }\big\}\]
(again, it is a club of $\lambda$). Finally, let 
\[C=\big\{\delta\in C^*\cap\mathop{\triangle}\limits_{\xi<\lambda} C_\xi: 
(\forall\xi,\zeta<\delta)(\vare(\{\xi,\zeta\})<\delta)\big\}.\]
Plainly, $C$ is a club of $\lambda$. Now, suppose that $\delta<\gamma$ are
two successive members of $C$. Put $Z_\delta=[\delta,\gamma)$ and let 
\[I_\delta=\big\{A\subseteq Z_\delta:\big(\exists\xi<\delta\big)\big(
\forall\alpha\in C^{p_\xi}\setminus\delta\big)\big(A\cap Z^{p_\xi}_\alpha
\notin d^{p_\xi}_\alpha\big)\big\}.\] 
It easily follows from the definition of $C$ that $I_\delta$ is a proper
ideal on $Z_\delta$, so we may pick an ultrafilter $d_\delta$ on $Z_\delta$
disjoint from $I_\delta$. Let $q=\big(C,\langle Z_\delta:\delta\in C\rangle,
\langle d_\delta:\delta\in C\rangle\big)$. Clearly $q\in\bqz$ and we will
argue that $q$ is a $\leq^0$--upper bound to $G^*$. So let $\xi<\lambda$. 
Suppose that $\delta\in C\setminus (\xi+1)$ and $A\in d_\delta$. Then
$A\notin I_\delta$, so there is $\alpha\in C^{p_\xi} \setminus\delta $ such
that $A\cap Z^{p_\xi}_\alpha\in d^{p_\xi}_\alpha$. Now we may use
\ref{tfcae} to conclude that $p_\xi\leq^0 q$.  
\end{proof}

\begin{proposition}
\label{1.4}
If $2^\lambda=\lambda^+$, then there is a $\leq^*_\bqz$--increasing sequence 
$\bar{p}=\langle p_\vare:\vare<\lambda^+\rangle\subseteq\bqz$ such that 
\[\fil(\bar{p})\stackrel{\rm def}{=}\bigcup\{\fil(p_\vare):\vare<
\lambda^+\}\]
is a uniform ultrafilter on $\lambda$.
\end{proposition}

\begin{proof}  
Straightforward induction using \ref{1.2}(4) and \ref{2.10}(2).
\end{proof}

For basic information on the ideal of meager subsets of ${}^\lambda\lambda$
and its covering number we refer the reader e.g. to Matet, Ros{\l}anowski
and Shelah \cite[\S 4]{MRSh:799}. Here we state only the definitions we will
need.  

\begin{definition}
\label{d0.1}
\begin{enumerate}
\item The space ${}^\lambda\lambda$ is endowed with the topology obtained by
taking as basic open sets $\emptyset$ and $O_s$ for $s\in{}^{\lambda>}
\lambda$, where $O_s=\{f\in {}^\lambda\lambda: s\subseteq f\}$. 
\item The $({<}\lambda)$--complete ideal of subsets of ${}^\lambda\lambda$
  generated by nowhere dense subsets of ${}^\lambda\lambda$ is denoted by
  $\Mkkk$.
\item ${\rm cov}(\Mkkk)$ is the minimal size of a family $\cA\subseteq\Mkkk$
  such that $\bigcup\cA={}^\lambda\lambda$.
\end{enumerate}
\end{definition}

\begin{theorem}
\label{covreason}
Assume that $\lambda^{<\lambda}=\lambda$ and ${\rm cov}(\Mkkk)=
2^\lambda$. Then there exists a very reasonable ultrafilter on $\lambda$.  
\end{theorem}

\begin{proof}
Fix a model $N\prec {\mathcal H}(\chi)$ (for some large regular cardinal
$\chi$) such that $\|N\|=\lambda$ and ${}^{\lambda>} N\subseteq N$.

For $p\in\bqz$ let $\langle\delta^p_\alpha:\alpha<\lambda\rangle$ be the
increasing enumeration of $C^p$ and let $\eta^p$ be the sequence of length
$\lambda$ such that 
\[\big(\forall\alpha<\lambda\big)\big(\eta^p(\alpha)=\langle
Z^p_{\delta^p_\alpha}, d^p_{\delta^p_\alpha}\rangle\big).\]
Next let 
\[\cT_\alpha=\big\{\eta^p\rest\alpha:p\in \bqz\big\}\cap N\quad\mbox{ (for
$\alpha<\lambda$)\quad and }\quad\cT=\bigcup\limits_{\alpha<\lambda}
\cT_\alpha.\] 
Clearly $\cT$ is a tree isomorphic to ${}^{\lambda>}\lambda$ by an
isomorphism preserving the levels (i.e., mapping $\cT_\alpha$ onto
${}^\alpha\lambda$). Also, every $\lambda$--branch $\eta\in\lim(\cT)$
determines a condition $p\in\bqz$ such that $\eta=\eta^p$. Let $Q^*=\{
p\in\bqz:\eta^p\in\lim(\cT)\}$. 

A family $G^*\subseteq Q^*$ is {\em linked\/} if it is $({<}\omega)$--linked
with respect to the partial order $\leq^0$<restricted to $Q^*$, that is if
every finite subset of $G^*$ has a $\leq^0$--upper bound in $Q^*$ (but the
bound does not have to be in $G^*$). Note that if $p_0,\ldots,p_n\in Q^*$
have a $\leq^0$--upper bound in $\bqz$, then they have a $\leq^0$--upper
bound in $Q^*$ as well.

For $p_0,\ldots,p_n\in G^*$, $\delta<\delta'<\lambda$ and an ultrafilter $d$
on $[\delta,\delta')$ let $(\oplus)^{p_0,\ldots,p_n}(\delta,\delta',d)$ mean 
\begin{enumerate}
\item[$(\oplus)^{p_0,\ldots,p_n}$]
\begin{enumerate}
\item $\delta,\delta'\in C^{p_0}\cap\ldots\cap C^{p_n}$, and 
\item if $B\in d$, $i\leq n$, then there is $\xi\in [\delta,\delta')\cap
  C^{p_i}$ such that $B\cap Z^{p_i}_\xi\in d^{p_i}_\xi$.
\end{enumerate}
\end{enumerate}
 
\begin{claim}
\label{cl1}
If $G^*\subseteq Q^*$ is a linked family, $\|G^*\|<{\rm cov}(\Mkkk)$, and 
$A\subseteq\lambda$, then there is $p\in Q^*$ such that 
\begin{enumerate}
\item[(a)] $G^*\cup\{p\}$ is linked, and 
\item[(b)] either $A\in \fil(p)$ or $\lambda\setminus A\in\fil(p)$.
\end{enumerate}
\end{claim}

\begin{proof}[Proof of the Claim]
We will consider two cases.
\smallskip 

\noindent{\sc Case 1:}\qquad {\em For every $p_0,\ldots,p_n\in G^*$,
  $n<\omega$, there is $p\in Q^*$ such that $A\in\fil(p)$ and $p_0\leq^0 
  p,\ldots,p_n \leq^0 p$.}
\smallskip 

Note that the assumption of the present case is equivalent to 
\begin{enumerate}
\item[$(\otimes)$] for every $p_0,\ldots, p_n\in G^*$, $n<\omega$, and
$\alpha <\lambda$ there are $\delta<\delta'<\lambda$ and an ultrafilter
$d\in N$ on $[\delta,\delta')$ such that $(\oplus)^{p_0,\ldots,p_n}(\delta,
\delta',d)$ holds true and $\alpha<\delta$ and $A\cap [\delta,\delta')\in
  d$.  
\end{enumerate}
(Remember that bounded subsets of $\lambda$ are in $N$.) We let 
\[\cT_A=\big\{\eta\in \cT:\big(\forall\alpha<\lh(\eta)\big)\big(\forall Z,d
\big)\big(\eta(\alpha)=\langle Z,d\rangle\ \Rightarrow\ A\cap Z\in d\big)
\big\}.\]
Clearly, $\cT_A$ is a $\lambda$--branching subtree of $\cT$ and $\cT_A$ is
isomorphic to ${}^{\lambda>}\lambda$. Now, for $p_0,\ldots,p_n\in G^*$,
$n<\omega$, and $\alpha<\lambda$ let  
\[\begin{array}{l}
I^A_\alpha(p_0,\ldots,p_n)\stackrel{\rm def}{=}\\
\qquad\big\{\eta\in \lim(\cT_A):\big(\exists \beta>\alpha\big)\big(\exists
\delta,\delta',d\big)\big((\oplus)^{p_0,\ldots,p_n}(\delta,\delta',d)\
\&\ \eta(\beta)=\langle[\delta,\delta'),d\rangle\big)\big\}.
  \end{array}\]
It should be clear that $I^A_\alpha(p_0,\ldots,p_n)$ is an open dense subset
of $\lim(\cT_A)$ (remember $(\otimes)$). Therefore (as $\|G^*\|<{\rm
cov}(\Mkkk)$) we know that 
\[\bigcap\big\{I^A_\alpha(p_0,\ldots,p_n):n<\omega\ \&\ p_0,\ldots,p_n\in
G^*\ \&\ \alpha<\lambda\big\}\neq\emptyset\]   
and we may choose $\eta$ from the set on the left-hand side above. Let $p\in 
Q^*$ be such that $\eta=\eta^p$. Since $\eta\in\lim(\cT_A)$ we know that
$A\in\fil(p)$. Also, for every $p_0,\ldots, p_n\in G^*$ we have $\eta\in
\bigcap\limits_{\alpha<\lambda} I^A_\alpha(p_0,\ldots,p_n)$ and hence 
\[\big\|\big\{\delta\in C^p:\mbox{ if }\delta'=\min\big(C^p\setminus
(\delta+1)\big)\mbox{ then }(\oplus)^{p_0,\ldots,p_n}(\delta,\delta',
d^p)\big\}\big\|=\lambda.\]
So one may easily construct $p^*\in Q^*$ which is $\leq^0$--stronger than 
$p,p_0,\ldots,p_n$ (remember \ref{tfcae}). Thus we have justified that
$G^*\cup\{p\}$ is linked.
\medskip

\noindent{\sc Case 2:}\qquad {\em There are $p_0,\ldots,p_n\in G^*$,
  $n<\omega$, such that
\begin{center}
if $p\in Q^*$ is $\leq^0$--stronger than $p_0,\ldots,p_n$, then
  $A\notin\fil(p)$.
\end{center}
} 
\smallskip 

It follows from the proof of \ref{1.2}(4) that then 
\begin{quotation}
{\em for every $q_0,\ldots,q_m\in G^*$, $m<\omega$, there is $q\in\bqz$ such
  that $\lambda\setminus A\in\fil(q)$ and $q_0\leq^0 q,\ldots, q_m\leq^0 q$}
\end{quotation}
(remember $G^*$ is linked). Thus we may repeat the arguments of Case 1 for
$\lambda\setminus A$ and find $p\in Q^*$ such that $G^*\cup\{p\}$ is linked 
and $\lambda\setminus A\in\fil(p)$. 
\end{proof}

\begin{claim}
\label{cl2}
If $G^*\subseteq Q^*$ is linked, $\|G^*\|<{\rm cov}(\Mkkk)$ and $p_0,p_1\in 
G^*$, then there is $p\in Q^*$ such that   
\begin{enumerate}
\item[(a)] $G^*\cup\{p\}$ is linked, and 
\item[(b)] $p_0\leq^0 p$ and $p_1\leq^0 p$. 
\end{enumerate}
\end{claim}

\begin{proof}[Proof of the Claim]
Let $p_0,p_1\in G^*$. Note that 
\begin{enumerate}
\item[$(\odot)$] for every $p_2,\ldots, p_n\in G^*$, $2\leq n<\omega$, and
$\alpha <\lambda$ there are $\delta<\delta'<\lambda$ and an ultrafilter
$d\in N$ on $[\delta,\delta')$ such that $(\oplus)^{p_0,p_1,p_2,\ldots,
p_n}(\delta,\delta',d)$ holds true and $\alpha<\delta$. 
\end{enumerate}
We let 
\[\cT^{p_0,p_1}=\big\{\eta\in \cT:\big(\forall\alpha<\lh(\eta)\big)\big(
\forall\delta,\delta',d\big)\big(\eta(\alpha)=\langle [\delta,\delta'),d
\rangle\ \Rightarrow\ (\oplus)^{p_0,p_1}(\delta,\delta',d)\big)\big\}\]
and we note that $\cT^{p_0,p_1}$ is a $\lambda$--branching subtree of $\cT$
isomorphic to ${}^{\lambda>}\lambda$. For $p_2,\ldots,p_n\in G^*$,
$2\leq n<\omega$, and $\alpha<\lambda$ we let  
\[\begin{array}{l}
I^{p_0,p_1}_\alpha(p_2,\ldots,p_n)\stackrel{\rm def}{=}\\
\big\{\eta\in \lim(\cT^{p_0,p_1}):\big(\exists \beta>\alpha\big)\big(
\exists\delta,\delta',d\big)\big((\oplus)^{p_2,\ldots,p_n}(\delta,\delta',d)\ 
\&\ \eta(\beta)=\langle[\delta,\delta'),d\rangle\big)\big\}.
  \end{array}\]
Then $I^{p_0,p_1}_\alpha(p_2,\ldots,p_n)$ is an open dense subset
of $\lim(\cT^{p_0,p_1})$ (remember $(\odot)$). Since $\|G^*\|<{\rm
cov}(\Mkkk)$, we may choose $p\in Q^*$ such that 
\[\eta^p\in\bigcap\big\{I^{p_0,p_1}_\alpha(p_2,\ldots,p_n):2\leq n<\omega\ \&\
p_2,\ldots,p_n \in G^*\ \&\ \alpha<\lambda\big\}\neq\emptyset.\]   
Like in the proof of \ref{cl1} we argue that $G^*\cup\{p\}$ is linked. Since
$\eta^p\in\lim(\cT^{p_0,p_1})$ we easily see that $p$ is $\leq^0$--stronger
than both $p_0$ and $p_1$.  
\end{proof}

\begin{claim}
\label{cl3}
If $G^*\subseteq Q^*$ is a linked family, $\|G^*\|<{\rm cov}(\Mkkk)$,
$\xi\leq\lambda$ is a limit ordinal and a sequence $\langle
p_\zeta:\zeta<\xi\rangle\subseteq G^*$ is $\leq^0$--increasing, then there
is $p\in Q^*$ such that   
\begin{enumerate}
\item[(a)] $G^*\cup\{p\}$ is linked, and 
\item[(b)] $(\forall\zeta<\xi)(p_\zeta\leq^0 p)$.
\end{enumerate}
\end{claim}

\begin{proof}[Proof of the Claim]
First let us consider the case when $\xi<\lambda$. Suppose that a sequence
$\bar{p}=\langle p_\zeta:\zeta<\xi\rangle\subseteq G^*$ is
$\leq^0$--increasing and let  
\[\cT_{\bar{p}}=\big\{\eta\in \cT:\big(\forall\alpha{<}\lh(\eta)\big)\big(
\forall \delta,\delta',d\big)\big(\eta(\alpha)=\langle [\delta,\delta'),
  d\rangle\ \Rightarrow\ (\forall\zeta{<}\xi)(\oplus)^{p_\zeta}(\delta,
\delta',d)\big)\big\}.\]
By arguments similar to that of \ref{1.2}(3) we verify that $\cT_{\bar{p}}$
is a $\lambda$--branching subtree of $\cT$ and it is isomorphic to
${}^{\lambda>}\lambda$. Like in the previous claims, for $p_0',\ldots,p_n'
\in G^*$, $n<\omega$, and $\alpha<\lambda$ we put 
\[\begin{array}{l}
I^{\bar{p}}_\alpha(p_0',\ldots,p_n')\stackrel{\rm def}{=}\\
\qquad\big\{\eta\in \lim(\cT_{\bar{p}}):\big(\exists\beta>\alpha\big)\big(
\exists\delta,\delta',d\big)\big((\oplus)^{p_0',\ldots,p_n'}(\delta,\delta',
d)\ \&\ \eta(\beta)=\langle[\delta,\delta'),d\rangle\big)\big\}.
  \end{array}\]
Then each $I^{\bar{p}}_\alpha(p_0',\ldots,p_n')$ is an open dense subset of 
$\lim(\cT_{\bar{p}})$. [Why? Let $\eta\in\cT_{\bar{p}}$. We may assume that
for each $\vare<\zeta<\xi$ and $\beta\in C^{p_\zeta}\setminus\lh(\eta)$ and
$A\in d^{p_\zeta}_\beta$ there is $\gamma\in C^{p_\vare}$ such that $A\cap 
Z^{p_\vare}_\gamma\in d^{p_\vare}_\gamma$. We also may demand that 
\[\delta_0\stackrel{\rm def}{=}\sup\big(\delta'<\lambda:\big(\exists\alpha
<\lh(\eta)\big)\big(\exists \delta,d\big)\big(\eta(\alpha)=\langle [\delta,
\delta'),d\rangle\big)\in \bigcap\limits_{\vare<\xi}C^{p_\vare}\cap
\bigcap\limits_{i\leq n} C^{p_i'}.\]  
Choose inductively a sequence $\langle\delta_\zeta,d_\zeta:\zeta<\xi
\rangle$ so that   
\begin{enumerate}
\item[(a)] $\langle\delta_\zeta:\zeta<\xi\rangle$ is an increasing
  continuous sequence of ordinals below $\lambda$, 
\item[(b)] $d_\zeta\in N$ is an ultrafilter on $[\delta_\zeta,\delta_{\zeta
    +1})$, $d_{\zeta+1}\in\bigcap\limits_{\vare<\xi} C^{p_\vare}$, 
\item[(c)] $(\oplus)^{p_0',\ldots,p_n',p_\zeta}(\delta_\zeta,\delta_{\zeta
  +1},d_\zeta)$ holds true (for each $\zeta<\xi$). 
\end{enumerate}
Let $\delta_\xi=\sup(\delta_\zeta:\zeta<\xi)$ and let $e\in N$ be an
ultrafilter on $\xi$ disjoint from the ideal of bounded subsets of
$\xi$. Put $d=\bigoplus\limits_{\zeta<\xi}^e d_\zeta$ --- it is an
ultrafilter on $[\delta_0,\delta_\xi)$, $d\in N$ and $(\oplus)^{p_0',
\ldots,p_n',p_\zeta}(\delta_0,\delta_\xi,d)$ holds true for each
  $\zeta<\xi$. Consequently $\eta\cup\big\{(\lh(\eta),\langle [\delta_0,
\delta_\xi),d\rangle)\big\}\in \cT_{\bar{p}}$ and every member of
$\lim(\cT_{\bar{p}})$ going through it belongs to $I^{\bar{p}}_\alpha(p_0,
\ldots,p_n)$.]    

Thus we may pick $p\in Q^*$ such that 
\[\eta^p\in \bigcap\big\{I^{\bar{p}}_\alpha(p_0,\ldots,p_n):n<\omega\ \&\
p_0,\ldots,p_n\in G^*\ \&\ \alpha<\lambda\big\}.\]   
Since $\eta^p\in\lim(\cT_{\bar{p}})$ we easily see that
$p_\zeta\leq^0 p$ for all $\zeta<\xi$, and like in the proof of \ref{cl1} we
argue that $G^*\cup\{p\}$ is linked. 
\smallskip 

If $\xi=\lambda$ and  $\bar{p}=\langle p_\zeta:\zeta<\lambda\rangle$ 
is $\leq^0$--increasing, then we proceed in a similar manner except that we
work in the tree
\[\cT^*_{\bar{p}}=\big\{\eta\in \cT:\big(\forall\alpha{<}\lh(\eta)\big)
\big(\forall \delta,\delta',d\big)\big(\eta(\alpha)=\langle [\delta,
\delta'),d\rangle\ \Rightarrow\ (\forall\zeta{<}\alpha)(\oplus)^{p_\zeta}(
\delta,\delta',d)\big)\big\}.\]
\end{proof}

\begin{claim}
\label{cl4}
Assume that $G^*\subseteq Q^*$ is a linked family, $\|G^*\|<{\rm
cov}(\Mkkk)$, $C\subseteq\lambda$ is a club and $\langle\delta_\xi:\xi<
\lambda\rangle$ is the increasing enumeration of $C$. Then there is $p\in
Q^*$ and a club $C^*$ of $\lambda$ such that  
\begin{enumerate}
\item[(a)] $G^*\cup\{p\}$ is linked, and 
\item[(b)] $\bigcup\big\{[\delta_{\xi+1},\delta_\zeta):\xi<\zeta$ are
  successive members of $ C^*\big\}\in \fil(p)$.  
\end{enumerate}
\end{claim}

\begin{proof}[Proof of the Claim]
Let  
\[\begin{array}{r}
\cT_C=\big\{\eta\in \cT:\mbox{for each }\alpha<\lh(\eta)\mbox{ such that
}\alpha=\delta_\alpha\mbox{ and for every }\alpha',d\qquad\\
\eta(\alpha)=\langle [\alpha,\alpha'),d\rangle\quad \Rightarrow\quad
  \delta_{\alpha+1}<\alpha'\ \&\ [\delta_\alpha,\delta_{\alpha+1})\notin
    d\big)\big\}. 
  \end{array}\]
One easily verifies that $\cT_C$ is a $\lambda$--branching subtree of $\cT$
which is isomorphic to ${}^{\lambda>}\lambda$. Like before, for $p_0,\ldots, 
p_n\in G^*$, $n<\omega$, and $\alpha<\lambda$ we put  
\[\begin{array}{l}
I^C_\alpha(p_0,\ldots,p_n)\stackrel{\rm def}{=}\\
\qquad\big\{\eta\in \lim(\cT_C):\big(\exists\beta>\alpha\big)\big(
\exists\delta,\delta',d\big)\big((\oplus)^{p_0,\ldots,p_n}(\delta,\delta',d)\
\&\ \eta(\beta)=\langle[\delta,\delta'),d\rangle\big)\big\}.
  \end{array}\]
Each $I^C_\alpha(p_0,\ldots,p_n)$ is an open dense subset of $\lim(\cT_C)$
and hence there is $p\in Q^*$ such that  
\[\eta^p\in \bigcap\big\{I^C_\alpha(p_0,\ldots,p_n):n<\omega\ \&\
p_0,\ldots,p_n\in G^*\ \&\ \alpha<\lambda\big\}.\]   
Like in the proof of \ref{cl1} we argue that $G^*\cup\{p\}$ is linked. Put 
\[C^*=\big\{\alpha<\lambda:\alpha=\delta_\alpha\mbox{ is limit } \&\
\big(\exists\alpha',d\big)\big(\eta^p(\alpha)=\langle [\alpha,\alpha'),d
\rangle\big)\big\}\]   
and note that $C^*$ is a club of $\lambda$. Note that if $\alpha\in C^*$ and 
$\eta^p(\alpha)=\langle [\alpha,\alpha'),d\rangle$, then
  $\delta_{\alpha+1}<\alpha'$ and $[\delta_\alpha,
\delta_{\alpha+1})\notin d$. Consequently, 
\[\bigcup\big\{[\delta_{\alpha+1},\delta_\beta):\alpha<\beta\mbox{ are
    successive members of }C^*\big\}\in \fil(p).\] 
\end{proof}

To prove the theorem we construct inductively a sequence $\langle q_\zeta:
\zeta<2^\lambda\rangle$ of elements of $Q^*$ such that  
\begin{itemize}
\item for each $\xi<2^\lambda$ the family $\{q_\zeta:\zeta<\xi\}$ is linked, 
\item for each $A\subseteq\lambda$ there is $\zeta<2^\lambda$ such that
  either $A\in\fil(q_\zeta)$ or $\lambda\setminus A\in\fil(q_\zeta)$,  
\item for each $\zeta<\xi<2^\lambda$ there is $\alpha<2^\lambda$ such that
  $q_\zeta\leq^0 q_\alpha$ and $q_\xi\leq^0 q_\alpha$,  
\item if $\xi\leq\lambda$ and $\langle p_\zeta:\zeta<\xi\rangle$ is a
  $\leq^0$--increasing sequence of elements of $\{q_\zeta:\zeta<2^\lambda\}$,
  then there is $\alpha<2^\lambda$ such that $q_\alpha$ is a $\leq^0$--upper
  bound to all $p_\zeta$'s, 
\item  if a sequence $\langle\delta_\xi:\xi<\lambda\rangle\subseteq\lambda$
  is increasing continuous, then for some $\zeta< 2^\lambda$ and a club
  $C^*$ of $\lambda$ we have 
\[\bigcup\big\{[\delta_{\xi+1},\delta_\zeta):\xi<\zeta\mbox{ are successive
  members of } C^*\big\}\in \fil(q_\zeta).\]  
\end{itemize}
The construction is a straightforward application of a suitable bookkeeping
device and Claims \ref{cl1}--\ref{cl4}. After it is carried out put
$G^*=\{q_\zeta:\zeta<2^\lambda\}$ and note that $\fil(G^*)$ is a very
reasonable ultrafilter on $\lambda$.  
\end{proof}

Let us finish this section with an observation showing that the assumption
$\lambda^{<\lambda}=\lambda$ in Theorem \ref{covreason} is very natural in
the given context.

\begin{proposition}
Assume $\theta<\lambda={\rm cf}(\lambda)<2^\theta$. Then ${\rm cov}(\Mkkk)
=\lambda^+$.
\end{proposition}

\begin{proof}
Let $\langle\nu_\xi:\xi<\lambda^+\rangle$ be a sequence of distinct
functions from $\theta$ to $2$. Let $\langle\delta_\alpha:\alpha<\lambda
\rangle\subseteq\lambda$ be an increasing continuous sequence such that
$\delta_0=0$, $\delta_{\alpha+1}=\delta_\alpha+\theta$ (for
$\alpha<\lambda$). Now, for $\xi<\lambda^+$ we define
\[F_\xi=\big\{\eta\in{}^\lambda\lambda:\big(\forall\alpha<\lambda\big)
\big(\exists i<\theta\big)\big(\eta(\delta_\alpha+i)\neq \nu_\xi(i)\big)
\big\}.\]
Plainly, each $F_\xi$ is a closed nowhere dense subset of
${}^\lambda\lambda$. We claim that $\bigcup\limits_{\xi<\lambda^+} F_\xi=
{}^\lambda\lambda$. To this end suppose that $\eta\in {}^\lambda\lambda$ and
consider the restrictions $\eta\rest [\delta_\alpha,\delta_{\alpha+1})$ for
$\alpha<\lambda$. These restrictions determine $\lambda$ functions from
$\theta$ to $2$, so we may find $\xi<\lambda^+$ such that $\nu_\xi$ is
distinct from all these functions, i.e., $(\forall\alpha<\lambda)(\exists
i<\theta)(\eta(\delta_\alpha+i)\neq\nu_\xi(i))$. Then $\eta\in F_\xi$.
\end{proof}

\section{$\fil(G^*)$ and dominating families}
In this section we show that families $G^*\subseteq\bqz$ generating
reasonable ultrafilters cannot be too small.  

\begin{theorem}
\label{bad}
For $p\in\bqz$ let $f_p\in {}^\lambda\lambda$ be such that 
\[\big(\forall\alpha<\lambda\big)\big(f_p(\alpha)\in C^p\ \&\
\otp\big(C^p\cap f_p(\alpha)\big)=\omega\cdot\alpha+\omega\big).\] 
\begin{enumerate}
\item Suppose that $G^*_0\subseteq\bqz$ is $({<}\aleph_1)$--directed (with 
respect to $\leq^0$) and $\fil(G^*)$ is a weakly reasonable ultrafilter. 
Then $\cF_0=\{f_p:p\in G^*_0\}$ is a dominating family in ${}^\lambda
\lambda$. 
\item Suppose that $G^*_1\subseteq\bqz$ is directed (with respect to
$\leq^0$) and $\fil(G^*_1)$ is a weakly reasonable ultrafilter on 
$\lambda$. Then $\cF_1=\{f_p:p\in G^*_1\}$ is a club--dominating family in
  ${}^\lambda\lambda$.
\end{enumerate}
\end{theorem}

\begin{proof}
(1)\quad First note that if $p,q\in G^*_0$, $p\leq^0 q$, then for some
$\vare<\lambda$, if $\alpha<\beta<\gamma$ are successive members of
$C^q\setminus\vare$, then $(\alpha,\gamma)\cap C^p\neq\emptyset$. Thus
$p\leq^0 q$ implies that for all sufficiently large $\alpha<\lambda$ we
have $f_p(\alpha)\leq f_q(\alpha)$. Consequently the family $\cF_0$ is
$({<}\aleph_1)$--directed (with respect to $\leq^*$). 

Suppose towards contradiction that $\cF_0$ is not a dominating
family. Then we may choose an increasing continuous sequence
$\bar{\alpha}^0=\langle\alpha^0_\xi:\xi<\lambda\rangle$ such that 
\[\big(\forall p\in G^*_0\big)\big(\exists^\lambda\vare<\lambda\big)\big(
f_p(\alpha^0_\vare)<\alpha^0_{\vare+1}\big).\]
Now, by induction on $n<\omega$, choose increasing continuous sequences
$\bar{\alpha}^n=\langle\alpha^n_\xi:\xi<\lambda\rangle$ so that letting
$C_n=\{\alpha^n_\xi:\xi<\lambda\}$ we have  
\begin{enumerate}
\item[(i)]   $\bar{\alpha}^0$ is the one chosen earlier,
\item[(ii)]  $C_{n+1}\subseteq\big\{\alpha^n_\vare:\vare=\alpha^n_\vare
\mbox{ is a limit ordinal }\big\}$,  
\item[(iii)] $\bigcup\big\{[\alpha^n_\vare,\alpha^n_{\vare+1}):\vare\in
  C_{n+1}\big\} \notin\fil(G^*_0)$. 
\end{enumerate}
It should be clear that the construction of $\bar{\alpha}^n$'s is possible
(remember that $\fil(G^*_0)$ is a weakly reasonable ultrafilter; use
\ref{easyob}). Let $C_\omega=\bigcap\limits_{n<\omega} C_n$ and let
$\langle\alpha^\omega_\xi:\xi<\lambda\rangle$ be the increasing enumeration
of $C_\omega$. It follows from (ii) that for every $\vare<\lambda$ the
sequence $\langle\alpha^n_{\vare+1}:n<\omega\rangle$ is strictly increasing
and $\sup(\alpha^n_{\vare+1}:n<\omega)\in C_\omega$, and if $\xi=
\alpha^\omega_\vare$, then $\sup(\alpha^n_{\xi+1}:n<\omega)=
\alpha^\omega_{\vare+1}$. It follows from (iii) that for every $n<\omega$ 
\[A_n\stackrel{\rm def}{=}\bigcup\big\{[\alpha^\omega_\vare,\alpha^n_{\xi
+1}):\vare<\lambda\ \&\ \xi=\alpha^\omega_\vare\big\}\notin\fil(G^*_0).\] 

Fix $p\in G^*_0$ for a moment. By the choice of $\bar{\alpha}^0$ we know
that the set $\{\xi<\lambda:f_p(\alpha^0_\xi)<\alpha^0_{\xi+1}\}$ is
unbounded in $\lambda$, and hence also the set $\{\vare<\lambda:
f_p(\alpha^\omega_\vare)<\alpha^\omega_{\vare+1}\}$ is unbounded in
$\lambda$. Therefore for some $n<\omega$ we have 
\[\|\{\vare<\lambda:\xi=\alpha^\omega_\vare\ \Rightarrow\ f_p(
\alpha^\omega_\vare)<\alpha^n_{\xi+1}\}\|=\lambda;\]
let $n(p)$ be the first such $n<\omega$. 

Note that if $p\leq^0 q$ are from $G^*_0$, then $n(p)\leq n(q)$ (as
$f_p\leq^* f_q$). Consequently, since $G^*_0$ is $({<}\aleph_1)$--directed,
there in $n^*<\omega$ such that $(\forall p\in G^*_0)(n(p)\leq n^*)$. Look
at the set $A_{n^*}$: for every $p\in G^*_0$ there are $\lambda$ many
$\vare<\lambda$ such that $\alpha^\omega_\vare<f_p(\alpha^\omega_\vare)<
\alpha^{n^*}_{\xi+1}$, where $\xi=\alpha^\omega_\vare$, and so (by the
definition of $f_p$) we get $A_{n^*}\in \big(\fil(p)\big)^+$. Since
$\fil(G^*_0)$ is an ultrafilter we get an immediate contradiction with
$A_{n^*}\notin \fil(G^*_0)$.  
\medskip

\noindent (2)\quad Suppose towards contradiction that $\cF_1$ is not
club--dominating in ${}^\lambda\lambda$. Then we may find an increasing
function $h\in {}^\lambda\lambda$ such that 
\[\big(\forall p\in G^*_1\big)\big(\big\{\vare<\lambda:f_p(\vare)<h(
\vare)\big\}\mbox{ is stationary in }\lambda\big).\]
Pick an increasing continuous sequence $\langle\delta_\xi:\xi<\lambda
\rangle\subseteq\lambda$ such that $(\forall \xi<\lambda)(h(\delta_\xi)<
\delta_{\xi+1})$. Since $\fil(G^*_1)$ is weakly reasonable, we may use
\ref{easyob} to pick a club $C$ of $\lambda$ such that $C\subseteq\{\xi<
\lambda:\delta_\xi=\xi\mbox{ is a limit ordinal }\}$ and 
\[\bigcup\big\{[\delta_\xi,\delta_{\xi+1}):\xi\in
  C\big\}\notin\fil(G^*_1).\] 
Since $\fil(G^*_1)$ is an ultrafilter, for some $p\in G^*_1$ we have 
\[\lambda\setminus\bigcup\big\{[\delta_\xi,\delta_{\xi+1}):\xi\in
  C\big\}\in\fil(p).\] 
However, by the choice of $h$, the set $\{\xi<\lambda:\delta_\xi=\xi\in C\
\&\ f_p(\xi)<h(\xi)<\delta_{\xi+1}\}$ is stationary (so of size $\lambda$),
and we get an immediate contradiction with the definition of $f_p$.
\end{proof}

\begin{corollary}
\begin{enumerate}
\item If $G^*_0\subseteq\bqz$ is $({<}\aleph_1)$--directed (with
respect to $\leq^0$) and $\fil(G^*)$ is a weakly reasonable ultrafilter,
then $\|G^*_0\|\geq \gdl$.
\item If $G^*_1\subseteq\bqz$ is directed (with respect to $\leq^0$) and
$\fil(G^*_1)$ is a weakly reasonable ultrafilter on $\lambda$, 
then $\|G^*_1\|\geq \gcl$.
\end{enumerate}
\end{corollary}

\begin{proposition}
\label{meas}
Suppose that $G^*_0\subseteq\bqz$ is $({<}\lambda)$--directed (with respect
to $\leq_\bqz$) and $\fil(G^*)$ is an ultrafilter. For $p\in G^*_0$ let
$f_p\in {}^\lambda\lambda$ be defined as in \ref{bad}. If $\cF_0=\{f_p:
p\in G^*_0\}$ is not a dominating family in ${}^\lambda\lambda$, then
$\lambda$ is measurable. 
\end{proposition}

\begin{proof}
Similarly as in the proof of \ref{bad}(1), we note that $\cF_0$ is
$({<}\lambda)$--directed (with respect to $\leq^*$). Assume $\cF_0$ is not
dominating family. Then may choose an increasing continuous sequence
$\langle\alpha_\xi:\xi<\lambda\rangle$ such that 
\[\big(\forall p\in G^*_0\big)\big(\exists^\lambda\vare<\lambda\big)\big(
f_p(\alpha_\vare)<\alpha_{\vare+1}\big).\]
Let 
\[\cU=\Big\{A\subseteq\lambda:\big(\exists p\in G^*_0\big)\big(\exists\delta<
\lambda\big)\big(\forall\vare>\delta\big)\big(f_p(\alpha_\vare)<
\alpha_{\vare+1}\ \Rightarrow\ \vare\in A\big)\Big\}.\]
We are going to show that $\cU$ is a $\lambda$--complete uniform ultrafilter
on $\lambda$.  It should be clear that $\cU$ includes all co-bounded subsets
of $\lambda$ and that it is a $\lambda$--complete filter (remember that
$\cF_0$ is ${<}\lambda$--directed). To show that it is an ultrafilter
suppose that $A\subseteq\lambda$ and let 
\[B=\bigcup\big\{[\alpha_\vare,\alpha_{\vare+1}):\vare\in A\big\}\subseteq
  \lambda.\] 
Since $\fil(G^*_0)$ is an ultrafilter, then either $B\in \fil(G^*_0)$ or
$\lambda\setminus B\in \fil(G^*_0)$. Suppose that the former happens, so we
may choose $p\in G^*_0$ such that $B\in\fil(p)$. Then for some
$\delta<\lambda$ we have  
\[\big(\forall\beta\in C^p\setminus\delta\big)\big(B\cap Z^p_\beta\in
d^p_\beta\big).\]
Now, if $\vare>\delta$ and $f_p(\alpha_\vare)<\alpha_{\vare+1}$, then
$[\alpha_\vare,\alpha_{\vare+1})\cap B\neq\emptyset$ and thus
$[\alpha_\vare,\alpha_{\vare+1})\subseteq B$, so $\vare\in A$. Consequently
$A\in\cU$ (as witnessed by $p,\delta$). In the same manner one shows that if
$\lambda\setminus B\in \fil(G^*_0)$, then $\lambda\setminus A\in\cU$.
\end{proof}

\section{Open problems and further investigations}
It may well be that our forcing techniques for uncountable $\lambda$ are
still not strong enough to carry out the arguments parallel to the
consistency results for ultrafilters on $\omega$. However, we feel that the
recent progress in the theory of forcing iterated with uncountable supports
(as exemplified by \cite{Sh:587}, Ros{\l}anowski and Shelah \cite{RoSh:655},
\cite{RoSh:777}, \cite{RoSh:860} and Eisworth \cite{Ei03}) may prove to be
useful in developing iterated forcing for ``killing'' and/or ``preserving''
some subfamilies of the class of reasonable ultrafilters. In particular, 
in Ros{\l}anowski and Shelah \cite{RoSh:890} we continue the research of the
present paper and we introduce {\em super reasonable ultrafilters\/}
which are stronger than very reasonable ultrafilters. We show that for
inaccessible $\lambda$ it is consistent that there are such ultrafilters
determined by generating systems of size less than $2^\lambda$, and we also
prove a result on preserving them in $\lambda$--support iterations. We also
show that consistently there are no super reasonable ultrafilters. These
results may be interpreted as some progress towards generalizing (a), (b)
and (c) from the Introduction. However, several other natural problems
remain untouched. One of the main questions we are interested in is 

\begin{problem}
\label{0.1}
Let $\lambda$ be a regular uncountable cardinal. Is it provable in ZFC that
  there exist reasonable ultrafilters on $\lambda$? Very reasonable? (See
  \ref{1.3}(4,5).) 
\end{problem}

\begin{problem}
\label{0.1F}
Is it consistent that there exists a very reasonable ultrafilter $D$ on
$\lambda$ such that for every very reasonable ultrafilter $D'$ on
$\lambda$ for some function $f\in\cFl$ we have $D/f=D'/f$? 
\end{problem}

Since in the present paper we deal with dividing by $f \in\cFl$, and the
normal ultrafilters are fixed points for this operation, the natural
question is:   

\begin{problem}
\label{0.1B}
Is it consistent that for every $D\in\uuf$ there is $f\in\cFl$ such that 
either $D/f$ is normal or $D/f$ is reasonable (or even very reasonable)? 
\end{problem}

We may also re-interpret our aim as follows. 

\begin{definition}
\label{1.3D}
\begin{enumerate}
\item Let $\UE^*_{\lambda,\mu}$ be the family of all $({<}\mu)$--directed
(with respect to $\leq^0$) subsets $G^*$ of $\bqz$ such that $\fil(G^*)$ is
  a proper ultrafilter on $\lambda$;  
\item $\UF^*_{\lambda,\mu}=\big\{\fil(G^*):G^*\in\UE^*_{\lambda,\mu}\big\}$;
$\UF^*_\lambda=\UF^*_{\lambda,\lambda^+}$ and $\UF_\lambda=
\UF^*_{\lambda,\aleph_0}$. 
\end{enumerate}
\end{definition}

\begin{aim}
\label{1.3E}
Investigate $\UF^*_\lambda,\UF_\lambda$; in particular can we have: any two
of them have common quotients.
\end{aim}

We expect that the forcing theorems needed for further research will be
similar to \cite{Sh:587} and more so to \cite{RoSh:655, RoSh:777, RoSh:860},
in some respects, and for others to \cite{Sh:832}.



\def\germ{\frak} \def\scr{\cal} \ifx\documentclass\undefinedcs
  \def\bf{\fam\bffam\tenbf}\def\rm{\fam0\tenrm}\fi 
  \def\defaultdefine#1#2{\expandafter\ifx\csname#1\endcsname\relax
  \expandafter\def\csname#1\endcsname{#2}\fi} \defaultdefine{Bbb}{\bf}
  \defaultdefine{frak}{\bf} \defaultdefine{=}{\B} 
  \defaultdefine{mathfrak}{\frak} \defaultdefine{mathbb}{\bf}
  \defaultdefine{mathcal}{\cal}
  \defaultdefine{beth}{BETH}\defaultdefine{cal}{\bf} \def\bbfI{{\Bbb I}}
  \def\mbox{\hbox} \def\text{\hbox} \def\om{\omega} \def\Cal#1{{\bf #1}}
  \def\pcf{pcf} \defaultdefine{cf}{cf} \defaultdefine{reals}{{\Bbb R}}
  \defaultdefine{real}{{\Bbb R}} \def\restriction{{|}} \def\club{CLUB}
  \def\w{\omega} \def\exist{\exists} \def\se{{\germ se}} \def\bb{{\bf b}}
  \def\equivalence{\equiv} \let\lt< \let\gt>
  \def\implies{\Rightarrow}\def\mathfrak{\bf}

\end{document}